\documentclass[
11pt]{amsart}
\usepackage{multicol}
\usepackage{epic}
\usepackage{epsfig}
\usepackage[english]{babel}

\usepackage{amsfonts}
\usepackage{amsthm}
\usepackage{mathrsfs}
\usepackage{colonequals}
\usepackage{amssymb,epsf}
\usepackage{amssymb,color}
\definecolor{c20}{rgb}{0.,0.7,0.}
\definecolor{c30}{rgb}{0.,0.,1.}
\definecolor{c40}{rgb}{1,0.1,0.7}
\definecolor{c50}{rgb}{1,0,0}
\definecolor{c60}{rgb}{0,0.9,0.1}
\newtheorem{Theorem}{Theorem}[section]
\newtheorem{Proposition}[Theorem]{Proposition}

\newtheorem{Lemma}[Theorem]{Lemma}

\newtheorem{Ex}{Example}[section]
\newtheorem{Rem}{Remark}[section]

\newcommand{\eq}{\colonequals}

\newcommand{\bb}{\mathbb}

\def \R {{\bb R}}
\def \N {{\bb N}}

\def \E{{\bb E}}

\def \pr{{\bb P}}

\def \vb{\vspace{2mm}}

\def\rw{\rightarrow}
\def\IF{\infty}
\newcommand{\BQN}{\begin{eqnarray}}
\newcommand{\EQN}{\end{eqnarray}}
\newcommand{\BQNY}{\begin{eqnarray*}}
\newcommand{\EQNY}{\end{eqnarray*}}

\newcommand{\halmos}{\hfill $\Box$}

\newcounter{mylistcnt}
\renewcommand{\themylistcnt}{{\rm[{\roman{mylistcnt}}]}}

\begin{document}

\title{L\'evy-driven GPS queues with heavy-tailed input}

\author{K. D\c{e}bicki}
\address{Krzysztof D\c{e}bicki, Mathematical Institute, University of Wroc\l aw, pl. Grunwaldzki 2/4, 50-384 Wroc\l aw, Poland}
\email{Krzysztof.Debicki@math.uni.wroc.pl}
\author{P. Liu}
\address{Peng Liu, Mathematical Institute, University of Wroc\l aw, pl. Grunwaldzki 2/4, 50-384 Wroc\l aw, Poland
and Department of Actuarial Science, University of Lausanne, UNIL-Dorigny 1015 Lausanne, Switzerland}
\email{liupnankaimath@163.com}
\author{M. Mandjes}
\address{Michel Mandjes, Korteweg-de Vries Institute for Mathematics,
University of Amsterdam, Amsterdam, the Netherlands}
\email{M.R.H.Mandjes@uva.nl}
\author{I. Sierpi\'nska-Tu\l acz}
\address{Iwona Sierpi\'nska-Tu\l acz, Mathematical Institute, University of Wroc\l aw, pl. Grunwaldzki 2/4, 50-384 Wroc\l aw, Poland}
\email{i.sierpinska@gmail.com}

\date{\today}
\maketitle

\begin{abstract}
In this paper we derive
exact large-buffer asymptotics for a two-class Generalized Processor Sharing (GPS) model, under the
assumption that the input traffic streams generated by both classes correspond to heavy-tailed L\'evy processes.
Four scenarios need to be distinguished, which differ in terms of (i)~the level of heavy-tailedness of the
driving L\'evy processes  as well as  (ii)~the values of the corresponding mean rates relative to the GPS weights.
The derived results are illustrated by two important special cases, in which the queues'
inputs are modeled by
heavy-tailed compound Poisson processes  and by
$\alpha$-stable L\'evy motions.

\vb


\end{abstract}

\section{Introduction}\label{s.introduction}
In queueing resources that are shared by multiple traffic streams, smooth streams potentially experience poor performance
when they are mixed with less regular streams. Indeed, under a first-come-first-serve (FCFS) discipline, users that correspond to a
highly variable input process may negatively affect the quality-of-service of other users. This motivates the attention paid
to more sophisticated  queueing disciplines, in which firm (per-user) performance guarantees can be given. One such a policy is the
{\it generalized processor sharing} (GPS) discipline. In GPS all users classes are guaranteed a certain service rate, whereas
the residual capacity is distributed according to a given allocation rule. The earliest (packet-based) implementations of GPS,
usually referred to as {\it weighted fair queueing} (WFQ), date
back to the late  1980s \cite{DEM}.

In many real-life systems, input streams may exhibit rather extreme types of irregularities.
For instance in the domain of communication networks, measurement studies show that
traffic patterns are typically heavy-tailed, in that there is a relatively high likelihood of an extremely large
amount of traffic being generated over a short time interval. Under FCFS all streams
would perceive roughly  the same performance, which is essentially determined by the input
class with the heaviest tail. GPS can be considered as a viable way to remedy this complication,
by offering each class a guaranteed service rate.

\vb


In this paper we consider a two-class GPS system, in which the inputs are
L\'evy processes with heavy-tailed
marginals; in our context, `heavy-tailed' refers to the complementary distribution function having a regularly varying tail.
This class of L\'evy processes covers many practically relevant processes.
In the first place, it contains the class of compound Poisson processes, in which independent
and identically distributed (i.i.d.) regularly-varying jobs arrive according to a Poisson process. In the second place, it covers the
class of $\alpha$-stable L\'evy motion; this class is particularly relevant, as it appears as the limiting process
for random walk models with increments that have infinite variance \cite{TaS94,WHI}.


Our main findings are the exact asymptotics of the tail distributions of both queues.
More specially, with $Q_i$ denoting the stationary workload of the $i$-th queue, we
find explicit functions $f_i(\cdot)$ such that ${\mathbb P}(Q_i>u)/f_i(u)
\to 1$ as $u\to\infty$; we write ${\mathbb P}(Q_i>u)\sim f_i(u)$.
As it turns out, depending on the interplay between the heaviness of both inputs'
tail distributions and
the stability of the queues while working in isolation,
one can distinguish four scenarios, each of them leading to qualitatively
different asymptotics. The resulting asymptotics lend themselves to an
intuitive explanation, in that they reveal the most likely way that
the workload under consideration exceeds $u$, for $u$ large. The proofs rely on
combining bounds that were derived earlier for related queueing systems, as well as
a set of newly derived inequalities. Related results for settings that are special cases of ours
can be found in e.g.\ \cite{BBJ,Lelarge}, whereas in \cite{MvU} the focus is on GPS systems
with Gaussian inputs.

The paper is organized as follows. Notation, assumptions and preliminaries are
presented in Section 2. Then Section 3 states the main results, in terms of the
exact asymptotics for all four scenarios. These results are used in Section 4 to
give the corresponding expressions for the compound Poisson and $\alpha$-stable
cases. All proofs are given in Section~5.

\section{Notation and model description}\label{s.model}

In this paper we consider a queueing system that consists of two queues and one server.
Each queue, which has infinite storage capacity, is fed by an own traffic class;
the corresponding input processes are assumed  to be mutually independent.
The total service rate of the server is $c>0$. Class $i$ is  assigned
a guaranteed service rate $\phi_i c>0$ (or `weight'), where $\phi_1 c + \phi_2 c = c$.
This effectively means that if both classes are backlogged, then class $i$ is served
at rate $\phi_i c$, for $i=1,2$. If class $i$ has no backlog,  then the other
class obtains the excess service rate.

Throughout this paper, we intensively use the concept of {\it cumulative input processes.}
We define by $Z_i(s,t)$ the cumulative
input to queue $i$ in interval $(s,t]$, for $i=1,2$ and $s<t$.
We assume that
\[
Z_i(s,t)=Z_i(t)-Z_i(s), \:\:\: i=1,2,
\]
where $\{Z_1(t):t\in \R \}$ and $\{Z_2(t):t\in \R \}$ are mutually independent
L\'evy processes.

As pointed out in the introduction, we specially consider the situation in which the L\'evy input processes are heavy-tailed. In
more concrete terms, this means that in the sequel we impose the following assumptions:

\vb

{\bf A1}~~$\pr(Z_{i}(1) > u) \sim u^{-\alpha_i}L_i(u)$, with
$\alpha_i>1$ and $L_i(\cdot)$ slowly varying at $\infty$, for $i=1,2$;

{\bf A2}~~$\E[Z_{i}(1)] = \mu_{i}$, with $\mu_1+\mu_2=\mu<c.$

\vb

We let
$\{Q_i(t):t\ge0\}$ denote the stationary buffer content processes
for class $i$, for $i= 1,2$.
Observe that
condition {\bf A2} guarantees stability of the system, implying existence of the
stationary buffer content processes. To shorten the notation we
throughout write
\[Q_i\stackrel{\rm d}{:=} Q_i(0),\:\:\:i =1,2.\]
Notice that the system's stability does not rule out that one of the queues `is in overload'
(if it would operate in isolation, that is):
one could have that $\mu_i>\phi_i c$ for one of the queues.

We denote by $B_i(s,t)$, for $i=1,2$,
the amount of service obtained by the $i${-th}
class in time interval $(s,t]$. Then there is the obvious identity\begin{eqnarray}\label{e2}
Q_i(t) = Q_i(s)+Z_i(s,t)-B_i(s,t), \quad \forall s<t.
\end{eqnarray}

According to Reich's formula \cite{REICH} (see also \cite{DvU} in the context of GPS queues)
we have the following distributional representation for the stationary workloads:
\[Q_i \stackrel{\rm d}{=} \sup_{t\geq 0} \{ Z_i(-t,0) - C_i(-t,0)\},\]
where $C_i(s,t)$ is the amount of the service available to class $i$ in the interval $(s,t]$.
The relation $C_i(s,t) \geq \phi_i c\,(t-s)$ holds for all $s<t$.

Additionally, it is convenient to introduce, for $\lambda_i>\mu_i$ and $\lambda>\mu$,
\[
Q_i^{\lambda_i}(t) :=\sup_{s\geq t}\{ Z_i(-s,t) - \lambda_i (t+s)\},\:\:\:\:
Q^{\lambda}(t):=\sup_{s\geq t} \{Z_1(-s,t)+Z_2(-s,t) - \lambda (t+s)\}.\]
Observe that
$Q_i^{\lambda_i} :=Q_i^{\lambda_i} (0)$ is distributed as the stationary buffer content of queue $i$
working in isolation, if it were served at rate $\lambda_i$ all the time. Likewise,
$Q^\lambda :=Q^{\lambda} (0)$ corresponds to the total stationary buffer content of the system, if it were
served at rate~$\lambda$.

Since the queues interact  
symmetrically, we  focus on
just
$\pr(Q_1>u),$ for $u\to\infty.$

\section{Main results}\label{s.main}
In this section we present the main results of the paper.
We distinguish four scenarios, that differ in terms of (i)~the heaviness of the individual
input processes, and (ii)~the individual queues being underloaded or overloaded.
The proofs of all the results presented in this section are relegated to Section \ref{s.proofs}.

\subsection{Second queue in overload}\label{s.case1}

We first consider the scenario
that the second queue is unstable when working in isolation: $\mu_2>\phi_2 c$.
In this case, if the input process of the second queue
generates traffic at its mean rate (which does not correspond to a rare event), then
it will be using its full guaranteed service rate. This pattern would
leave the first queue as if working in isolation.
Based on this observation, we expect that
\[\pr(Q_{1} > u) \sim \pr(Q_{1}^{\phi_1c}>u).\]
The following theorem formalizes this heuristic. Notice that in this scenario we
necessarily have $\mu_1<\phi_1 c.$

\begin{Theorem}\label{th.secUnstab}
Suppose that $Z_1,Z_2$ satisfy {\bf A1-A2}.
If $\mu_2>\phi_2 c$, then, as $u\to\infty$,
$$\pr(Q_{1} > u) \sim \frac{1}{(\phi_1 c - \mu_1)(\alpha_1-1)}u^{1-\alpha_1}L_1(u).$$
\end{Theorem}

\subsection{Second queue in underload, first class is heavier}
\label{s.case2}
In the other three scenarios the second queue is stable while working in isolation, i.e.,
we consider the situation that
$\mu_2<\phi_2 c$.
As it turns out, under this condition the interplay between both input processes
plays a key role. We first concentrate on the case that
the first class is heavier than the second one, i.e., $\alpha_1<\alpha_2.$
Since the second queue is stable while working in isolation
and `is lighter' than the first one, the most likely way to generate a large workload in the
first queue does not involve a large  buffer
content in the second queue. The most probable way the first buffer
reaches a large level corresponds to (i)~the second class generating traffic
at its mean level $\mu_2$, and (ii)
the remaining service capacity $c-\mu_2$ being allocated to the first queue.
Hence the so-called \textit{reduced-load}
equivalence holds in this case, cf.\ e.g.\ \cite{BBJ2}:
\[\pr(Q_{1} > u) \sim \pr(Q_{1}^{c-\mu_2}>u).\]

This leads to the following theorem.

\begin{Theorem}\label{th.secStabH}
Suppose that $Z_1,Z_2$ satisfy {\bf A1-A2}.
If $\mu_2<\phi_2 c$ and $\alpha_1<\alpha_2$, then, as $u\to\infty$,
$$\pr(Q_{1} > u) \sim \frac{1}{(c-\mu)(\alpha_1-1)}u^{1-\alpha_1}L_1(u).$$
\end{Theorem}

\subsection{Second queue in underload, second class is heavier}
\label{s.case3}
In the remaining two scenarios the second queue is stable while working in isolation, and
the second class is heavier than the first one, i.e., $\alpha_2<\alpha_1.$ Two cases still need to be distinguished:
the first queue being in underload or not.

In this third scenario we suppose that both the first and the second queue are stable while working in isolation, i.e., $\mu_i<\phi_i c$ for $i=1,2$
(and, as mentioned above, the second class is the heavier).
For this scenario it turns out that again the reduced load equivalence holds:
\[\pr(Q_{1} > u) \sim \pr(Q_{1}^{c-\mu_2}>u).\]
Intuitively, this means that the most probable way in which queue 1 grows large
is as follows: the second class generates traffic roughly at its mean rate, and
the first queue builds up as acting in isolation
with service rate $c-\mu_2$ (which can be interpreted as the service rate left by the second queue).
Although the asymptotics coincide with those obtained in Theorem \ref{th.secStabH},
the proof of the upper bound for this case needs an entirely different approach (which motivates why we treat
them as separate cases).

\begin{Theorem}\label{th.secStabL}
Suppose that $Z_1,Z_2$ satisfy {\bf A1-A2}.
If $\mu_1<\phi_1 c$, $\mu_2<\phi_2 c$ and $\alpha_2<\alpha_1$, then, as $u\to\infty$,
$$\pr(Q_{1} > u) \sim  \frac{1}{(c-\mu)(\alpha_1-1)}u^{1-\alpha_1}L_1(u).$$
\end{Theorem}

\subsection{First queue in overload, second class is heavier}\label{s.case4}
Finally, we consider the scenario that
the first queue is in overload (i.e., unstable when working in isolation: $\mu_1>\phi_1 c$), and
the second class is the heavier (i.e., $\alpha_2<\alpha_1$).
We in addition assume that $Z_2$ be spectrally positive.

In this case
the most probable way in which the first queue reaches a high level is
such that the first class generates traffic roughly at its average rate $\mu_1$ (which does not correspond to a rare event).
Now the crucial issue concerns the fraction of its service rate that is left by the second class to the first class.
As it turns out, the most likely behavior of the second queue
can be linked with
the downstream queue of a fictitious two-node tandem queue
fed by $Z_2$ with service rate $\phi_2c$ at the upstream queue and
service rate $c-\mu_1$
at the downstream queue, in the sense that
\[\pr(Q_{1} > u) \sim \pr(V>u),\]
where
\[V \eq \sup_{t \geq 0}\{Z_2(-t,0) - (c-\mu_1)t\} - \sup_{s \geq 0} \{Z_2(-s,0) - \phi_2 cs\}.\]
This relation was also observed in
GPS models with fractional Brownian input in \cite{DMGPS}, whereas \cite{BBJ}
finds  a similar relation for the case of heavy tailed on-off input.
Combining this with results from \cite{DMvUTand} on tandem queues with spectrally positive input,
we thus arrive at the following asymptotics.

\begin{Theorem}\label{th.reqc}
Suppose that $Z_1,Z_2$ satisfy {\bf A1-A2},
$\mu_1>\phi_1 c$, $\alpha_2<\alpha_1$ and $Z_2$ is spectrally positive
with $\alpha_2\notin \N$.
Then, as $u\to\infty$,
\[\pr(Q_{1} > u) \sim \left( \frac{\mu_1-\phi_1c}{\phi_2c-\mu_2} \right)^{\alpha_2-1}\frac{1}{(c-\mu)(\alpha_2-1)}u^{1-\alpha_2}L_2(u).\]
\end{Theorem}

\begin{Rem}{\em
In the proof of
Theorem \ref{th.reqc} it plays a crucial role that $Z_2$ is assumed to be spectrally positive.
We strongly believe that this assumption is of a technical nature, in that  the statement of
Theorem
\ref{th.reqc} is valid for general $Z_2$.
We anticipate, however, that a proof for general $Z_2$ would be considerably more complicated,
and would go along entirely different lines;   see also Remark \ref{rem}.}
\end{Rem}

Observe that in the first three scenarios the workload of the first queue inherits the tail behavior of its input process: the complementary distribution function ${\mathbb P}(Q_1>u)$ essentially behaves as
$u^{1-\alpha_1}$. We conclude that in these cases the GPS mechanism succeeds in protecting the first stream against the second stream. Only in the last scenario,  ${\mathbb P}(Q_1>u)$ becomes heavier, which issue to the relatively large weight allocated to the second stream.

\section{Special cases}\label{s.examples}
In this section we use the general results, as presented in the previous section,
to find the asymptotics for
GPS systems fed by   compound Poisson processes with heavy-tailed
jumps (Section \ref{ss.comp})  and by $\alpha$-stable L\'evy input (Section \ref{ss.stable}).

\subsection{Compound Poisson input}\label{ss.comp}
This subsection concentrates on the case of {compound Poisson} inputs.
More concretely, we assume that $Z_i(t)$ is of the form
\[Z_i(t)=\sum_{k=1}^{N_i(t)}B_{k,i},\:\:\:\:
i=1,2.\]
In this definition of $Z_i(t)$, we assume that the processes $N_i(\cdot)$  are independent Poisson processes with rates $\lambda_i>0$.
In addition,
$(B_{k,1})_k$ and  $(B_{k,2})_k$ are both sequences of i.i.d.\ non-negative random variables, which are
independent of the processes $N_1(\cdot)$ and $N_2(\cdot)$.
We denote by $B_1, B_2$ the
generic random variables corresponding to the sequences $B_{k,1}$ and  $B_{k,2}$, where $F_1(\cdot)$ and $F_2(\cdot)$ denote their respective distribution functions.

The following proposition translate the
findings of the previous section into the setting of the
compound Poisson input model.

\begin{Proposition}
Assume that both $Z_1(\cdot)$ and $Z_2(\cdot)$ are independent compound Poisson processes
with $\mu_i=\E[Z_i(1)]=\lambda_i\E[B_i]$ and
$1-F_i(x)\sim x^{-\alpha_i}L_i(x)$, as $x\to\infty$, for $i=1,2$ and $L_i(\cdot)$
being slowly varying at $\infty$, with $\alpha_i>1$.

1) If $\mu_2>\phi_2 c$, then, as $u\to\infty$,
$$\pr(Q_1>u) \sim \frac{\lambda_1}{\phi_1c-\mu_1}\frac{1}{\alpha_1 -1}u^{1-\alpha_1}L_1(u).$$
\\
2) If $\mu_2<\phi_2 c$ and $\alpha_1<\alpha_2$, then, as $u\to\infty$, $$\pr(Q_1>u) \sim \frac{\lambda_1}{c-\mu}\frac{1}{\alpha_1 -1}u^{1-\alpha_1}L_1(u).$$
3) If $\mu_2<\phi_2 c$, $\mu_1<\phi_1 c$ and $\alpha_2<\alpha_1$,
then, as $u\to\infty$, $$\pr(Q_1>u) \sim \frac{\lambda_1}{c-\mu}\frac{1}{\alpha_1 -1}u^{1-\alpha_1}L_1(u).$$
4) If $\mu_1>\phi_1 c$, $\alpha_2<\alpha_1$ and $Z_2$ is spectrally positive, then, as $u\to\infty$, $$\pr(Q_1>u) \sim \frac{\lambda_2}{c-\mu}\left( \frac{\mu_1-\phi_1c}{\phi_2c-\mu_2}
\right)^{\alpha_2-1}\frac{1}{\alpha_2 -1}u^{1-\alpha_2}L_2(u).$$
\end{Proposition}

\proof
The proof follows straightforwardly from
Theorems \ref{th.secUnstab}, \ref{th.secStabH}, \ref{th.secStabL} and \ref{th.reqc}, respectively, in combination with Theorem 2.1 in \cite{ASM}.
\endproof

\subsection{$\alpha$-stable L\'evy input}\label{ss.stable}
In this second subsection we focus on the special case of $Z_1(\cdot)$ and $Z_2(\cdot)$
being independent  $\alpha_j$-stable L\'evy motions. This formally means that its law is given in terms of its characteristic function:
\[
\log\E e^{i\theta Z_j(1)} =
- |\theta|^{\alpha_j} (1-i
\beta_j{\rm sign}(\theta)\tan(\pi\alpha_j/2))+i \mu_j\theta,
\]
where $\alpha_j\in(1,2]$, $\beta_j\in(-1,1]$,
$\mu_j\in{\mathbb R}$, and ${\rm
sign}(x):=1_{(0,\infty)}(x)-1_{(-\infty,0)}(x).$
We write $Z_j\in \mathbb{S}(\alpha_j,\beta_j,\mu_j),$
see e.g., \cite{TaS94} or \cite{DeM15}.

Using that
\[\pr(Z_j(1)>x) \sim
c_{\alpha_j}(1+\beta_j)x^{-\alpha_j},\:\:\:\:\mbox{}\:\:\:\:c_{\alpha} :=
\frac{1-\alpha}{2\Gamma(2-\alpha)\cos (\pi\alpha/2)},\]
see \cite{TaS94}, in combination with the results presented in Section \ref{s.main},
we arrive at the following proposition.

\begin{Proposition}
\label{alphacol}
Suppose that $Z_i\in \mathbb{S}(\alpha_i,\beta_i,\mu_i),$  with $\alpha_i\in(1,2)$ for $i=1,2$.

1) If $\mu_2>\phi_2 c$, then, as $u\to\infty$, $$\pr(Q_1>u) \sim \frac{c_{\alpha_1}(1+\beta_1)}{(\phi_1c-\mu_1)(\alpha_1-1)}u^{1-\alpha_1}.$$
2) If $\mu_2<\phi_2 c$ and $\alpha_1<\alpha_2$, then, as $u\to\infty$, $$\pr(Q_1>u) \sim \frac{c_{\alpha_1}(1+\beta_1)}{(c-\mu)(\alpha_1-1)}u^{1-\alpha_1}.$$
3) If $\mu_2<\phi_2 c$, $\mu_1<\phi_1 c$ and $\alpha_2<\alpha_1$, then, as $u\to\infty$, $$\pr(Q_1>u) \sim \frac{c_{\alpha_1}(1+\beta_1)}{(c-\mu)(\alpha_1-1)}u^{1-\alpha_1}.$$
4) If $\mu_1>\phi_1 c$, $\alpha_2<\alpha_1$ and $\beta_2 = 1$, then, as $u\to\infty$,
$$\pr(Q_1>u) \sim \frac{2c_{\alpha_2}}{(c-\mu)(\alpha_2-1)}
\left( \frac{\mu_1-\phi_1c}{\phi_2c-\mu_2}
\right)^{\alpha_2-1}u^{1-\alpha_2}.$$
\end{Proposition}

\begin{Rem}\label{rem}{\em
Complementary to case 4) of Proposition \ref{alphacol},
for $\beta_2 \in (-1,1]$ (and $\mu_1>\phi_1 c$, $\alpha_2<\alpha_1$),
we can find asymptotic upper and lower bounds on $\pr(Q_1>u)$
that are tight up to a constant.
In particular,
combining the proof of Theorem \ref{th.reqc} with Theorem 5.3 and  Lemma 5.4 in \cite{DMvUTand} we obtain,
as $u\to\infty$,
$$\limsup_{u\to\infty}
\pr(Q_1>u)u^{-(1-\alpha_2)} \leq
\left( \frac{c_{\alpha_2}(1+\beta_2)}{(c-\mu)(\alpha_2-1)}+
\frac{c_{\alpha_2}(1+\beta_2)}{\phi_2c}\right)
\left( \frac{\mu_1-\phi_1c}{\phi_2c}
\right)^{\alpha_2-1},$$\\
and
$$\liminf_{u\to\infty}
\pr(Q_1>u)u^{-(1-\alpha_2)}\geq  \frac{c_{\alpha_2}(1+\beta_2)}{(c-\mu)(\alpha_2-1)}\left( \frac{\mu_1-\phi_1c}{\phi_2c}
\right)^{\alpha_2-1}.$$}
\end{Rem}


\section{Proofs}
\label{s.proofs}

Before we provide detailed proofs of the
results of Section \ref{s.main}, we present some
useful lemmas.
We begin with the classical result by Port \cite{Por89}, describing the asymptotics of the tail distribution of a single queue that is emptied at rate $c$.

\begin{Lemma}
\label{th.Port}
Suppose that $Z_1$ satisfies {\bf A1}, {\bf A2} with
$c>\mu_1$. Then, as $u\to\infty$,
$$\pr(Q_1^c > u) \sim \frac {1}{c-\mu_1}\frac{1}{\alpha_1-1}u^{1-\alpha_1}L_1(u)
.$$
\end{Lemma}

The following result is due to Willekens \cite{Wil87}, describing the asymptotic distribution of the
supremum of $Z_1(-t,0) - ct$ over a finite interval.

\begin{Lemma}
\label{th.Willekens}
Suppose that $Z_1$ satisfies {\bf A1}. Then, for each $T>0$, as $u\to\infty$,
$$
\pr \left( \sup_{t \in [0, T]}\{ Z_1(-t,0) - ct\}>u \right)\sim \pr (Z_1(1)>u).
$$
\end{Lemma}

Whereas the previous lemma considers the supremum over a finite interval, in the next lemma
the interval grows with the exceedance level $u$. This result may have appeared in some form in the literature, but we decided to include it here, as it has a natural and insightful proof.

\begin{Lemma}\label{lemma2}
Suppose that $Z_1$ satisfies {\bf A1}
with  $c>\mu_1$ and $\lim_{u\to\infty}{T(u)}/{u}=\infty$. Then, as $u\to\infty$,
$$\pr(Q_1^c > u)\sim \pr\left(\sup_{t\in [0, T(u)]}\{Z_1(t)-ct\}>u\right).$$
\end{Lemma}
\proof
Observe that the following trivial inequality holds:
\[
{\mathbb P}({\mathscr E}(u))\leq \pr(Q_1^c > u)\leq {\mathbb P}({\mathscr E}(u))+{\mathbb P}({\mathscr F}(u)),\]
with
\[{\mathscr E}(u):=\left\{
\sup_{t\in [0, T(u)]}\{Z_1(t)-ct\}>u\right\},\:\:\:\:
{\mathscr F}(u):=\left\{\sup_{t\geq T(u)}\{Z_1(t)-ct\}>u\right\}.
\]
Let $\widetilde{Q}_1^{c}\stackrel{\rm d}{=} Q_1^{c}$, with $\widetilde{Q}_1^{c}$ being independent of $\{Z_1(t), t\in\mathbb{R}\}$.
Then using the independence and stationarity of the increments of $Z_1$, we have, with $\varepsilon\in(0, c-\mu_1)$,
\BQNY
{\mathbb P}({\mathscr F}(u))
&=&\pr\left(Z_1(T(u))-cT(u)+\sup_{t\geq T(u)}\{Z_1(t)-Z_1(T(u))-c(t-T(u))\}>u\right)\\
&=&\pr\left(Z_1(T(u))+\widetilde{Q}_1^{c}>u+cT(u)\right)\\
&=&\pr\left(Z_1(T(u))-(\mu_1+\varepsilon)T(u)+\widetilde{Q}_1^{c}>u+(c-\mu_1-\varepsilon)T(u)\right).
\EQNY
By applying $\pr(X+Y\geq z) \leq \pr(X\geq fz) +\pr(Y\geq (1-f)z)$ for $f\in(0,1)$, this quantity is in turn bounded from above by, with $\Delta:=c-\mu_1-\varepsilon>0$,
\[\pr\left(Z_1(T(u))-(\mu_1+\varepsilon)T(u)\geq \frac{1}{2}\,\Delta T(u)\right)
+\pr\left({Q}_1^{c}>u+\frac{1}{2}\,\Delta T(u)\right).\]
We prove for each of these probabilities that they are $o\left(\pr(Q_1^c > u)\right)$
as $u\to\infty$.
The first probability is majorized by
\[\pr\left(\sup_{t\geq 0}\{Z_1(t)-(\mu_1+\varepsilon)t\}\geq \frac{1}{2}\,\Delta T(u)\right) =
\pr\left(Q_1^{\mu_1+\varepsilon} \geq \frac{1}{2}\,\Delta T(u)\right),\]
Now it follows from Lemma \ref{th.Port} that, recalling that $u=o(T(u))$,
$$\pr\left(Q_1^{\mu_1+\varepsilon} \geq \frac{1}{2}\,\Delta T(u)\right)=o\left(\pr(Q_1^c > u)\right), \ \ \pr\left({Q}_1^{c}>u+\frac{1}{2}\,\Delta T(u)\right)=o\left(\pr(Q_1^c > u)\right). $$
Therefore, we conclude that
${\mathbb P}({\mathscr F}(u))=o\left(\pr(Q_1^c > u)\right),$
which completes the proof. \endproof

Define, for $\lambda<\mu_2,$  ${T_\lambda}(u):=u/\bar {T_\lambda}(u)$, where
\[\bar T_{\lambda}(u):={\sqrt{\pr \left(\check Q_2^\lambda(0)>u/2\right)\vee (1/\log u) }},\:\:\:
\check Q_2^{\lambda}(s):=\sup_{t\geq s}\{Z_2(s)-Z_2(t)+\lambda (t-s)\}.\]

\begin{Lemma}\label{Lemma1} Suppose that $Z_2$ satisfies {\bf A1} with $\lambda<\mu_2$. Then,
as $u\to\infty$,
\[\xi(u):=\pr\left(\sup_{s\in [0,{T_\lambda}(u)]}\check Q_2^{\lambda}(s)>u\right)\to 0.\]
\end{Lemma}
\proof It is immediate  that, with $W(t):=Z_2(t)-\lambda t$,
\BQNY
\sup_{s\in [0,{T_\lambda}(u)]}\check Q_2^{\lambda}(s)&=&\sup_{s\in [0,{T_\lambda}(u)]}\sup_{t\geq s}\{Z_2(s)-Z_2(t)+\lambda (t-s)\}\\
&\leq& \sup_{s\in [0,{T_\lambda}(u)]}\sup_{t\geq 0}\{Z_2(s)-Z_2(t)+\lambda (t-s)\}=\check Q_2^{\lambda}(0)+\sup_{s\in[0,{T_\lambda}(u)]}W(s).
\EQNY
Setting $$g(u):=u\left(\bar {T_\lambda}(u) \right)^{1/2}, \ \
h(u):=u\left(\bar {T_\lambda}(u) \right)^{1/3},$$
we have
\BQNY
\zeta(u)&:=&\pr\left(\sup_{t\in [0,g(u)]}{\check Q}_2 ^{\lambda}(t)>u\right) \leq \pr\left({\check Q}_2 ^{\lambda}(0)+\sup_{s\in[0, g(u)]}W(s)>u\right)\\
&\leq& \pr\left({\check Q}_2 ^{\lambda}(0)+\sup_{s\in[0,g(u)]}W(s)>u, \sup_{s\in[0,g(u)]}W(s)\leq h(u)\right)\\
&& +\, \pr\left({\check Q}_2 ^{\lambda}(0)+\sup_{s\in[0,g(u)]}W(s)>u, \sup_{s\in[0,g(u)]}W(s)> h(u)\right)\\
&\leq&   \pr\left({\check Q}_2 ^{\lambda}(0)>u-h(u)\right)+\pr\left(\sup_{s\in[0,g(u)]}W(s)> h(u)\right)\\
&\leq&  \pr\left({\check Q}_2 ^{\lambda}(0)>u/2\right)+\pr\left(\sup_{s\in[0,g(u)]}\{Z_2(s)-(\mu_2+1)s\}> h(u)-(\mu_2-\lambda+1)g(u)\right)\\
&\leq&  \pr\left({\check Q}_2 ^{\lambda}(0)>u/2\right)+\pr\left(Q_2^{\mu_2+1}> h(u)/2\right).
\EQNY
Using the stationarity of ${\check Q}_2 ^{\lambda}(t)$, and noting that $m(u):=T_\lambda(u)/g(u)\to\IF$ as $u\to\IF$,
\[\xi(u) \leq \sum_{i=1}^{[ m(u)]+1}\pr\left(\sup_{t\in[(i-1)g(u),\,i\,g(u)]}{\check Q}_2 ^{\lambda}(t)>u\right)=([ m(u)]+1) \,\zeta(u).\]
Now, by applying  Lemma \ref{th.Port}, we have (noting that $h(u)\to\IF$ as $u\to \IF$)
that there is a positive constant $\kappa$ such that
\BQNY
\xi(u)&\leq& \frac{2{T_\lambda}(u)}{g(u)}\zeta(u)
\leq\frac{2{T_\lambda}(u)}{g(u)}\left(\pr\left({\check Q}_2 ^{\lambda}(0)>u/2\right)+\pr\left(Q_2^{\mu_2+1}> h(u)/2\right)\right)\\
&\leq& 2\left(\pr\left({\check Q}_2 ^{\lambda}(0)>u/2\right)\right)^{1/4}+\kappa\,(\log u)^{3/4}\,L_2(h(u))(h(u))^{1-\alpha_2}\to 0.
\EQNY
This completes the proof. \endproof

The following lemma plays an important role in the proof of Theorem \ref{th.secStabL}.
\begin{Lemma}\label{BIG}
Suppose that $Z_2$ satisfies {\bf A1} with  $\lambda<\mu_2$.  Then
$$\frac{Q_2^\lambda(t)}{t}\rw 0, \:\:\:{\rm a.s.},\:\:\: {\rm as}\:  t\rw\IF.$$
\end{Lemma}
\proof It suffices to prove that for some $\beta>0$,
\BQN\label{e1}
\sup_{t\in [n^\beta, (n+1)^\beta]}\frac{Q_2^{\lambda}(t)}{n^{\beta}}
\rw 0, \:\:\:{\rm a.s.},\:\:\: {\rm as}  \: n\rw\IF.
\EQN
Let  $\beta>({\alpha_2-1})^{-1}>0$ hereafter.
Then, with $n_+:=n+1$, $I_n:=[n^\beta, n_+^\beta]$,
\BQNY
\sup_{t\in I_n}Q_2^{\lambda}(t)&=&
\sup_{t\in I_n}\sup_{s\leq t}\{Z_2(t)-Z_2(s)-\lambda(t-s)\}\\
&\leq &\sup_{t\in I_n}\sup_{s\leq n_+^\beta}\{Z_2(n_+^\beta)-Z_2(s)-\lambda(n_+^\beta-s)+Z_2(t)-Z_2(n_+^\beta)  -\lambda(t-n_+^\beta)\}\\
&=&
Q_2^\lambda(n_+^\beta)+\sup_{t\in I_n}\{W(t)-W(n_+^\beta)\},
\EQNY
where, as before, $W(t):=Z_2(t)-\lambda t$.
In light of Lemma \ref{th.Port} we have that for any $\varepsilon>0$ and $\gamma\in(1, \beta(\alpha_2-1))$ there exists $N_0\in{\mathbb N}$   such that (where we recall that $\beta>({\alpha_2-1})^{-1}$)
\[\sum_{n=N_0}^\IF\pr\left(\frac{Q_2^\lambda(n_+^\beta)}{n^{\beta}}>\varepsilon\right)=\sum_{n=N_0}^\IF\pr\left(Q_2^\lambda(0)>\varepsilon n^{\beta}\right)\leq  \sum_{n=N_0}^\IF n^{-\gamma}<\IF.\]
Using the Borel-Cantelli lemma, we obtain that $Q_2^\lambda(n_+^\beta)/n^\beta \rw 0$ a.s.\ as $n\rw\IF.$
In order to establish (\ref{e1}) we are now left to prove that $n^{-\beta} \sup_{t\in I_n}\{W(t)-W(n_+^\beta)\} \rw 0$ a.s., as $n\rw\IF.$ This convergence is established as follows.

By the fact that
${Z_2(t)}/t  \rw \mu_2,$ a.s., as $t\rw\IF$, we have ${W(t)}/t  \rw \mu_2-\lambda,$ a.s., as $t\rw\IF$, and
$$\sup_{s,t\geq n}\left\{\frac{W(s)}{s}-\frac{W(t)}{t}\right\}\rw 0, \:\mbox{a.s.},\:\:\mbox{as $n\rw\IF$}.$$
Therefore,
\BQN
0&\leq&\frac{1}{n^{\beta}}\sup_{t\in I_n}\{W(t)-W(n_+^\beta) \}\leq\frac{n_+^\beta}{n^\beta}\sup_{t\in I_n}\left\{\frac{W(t)}{t}-\frac{W(n_+^\beta)}{t}\right\}
\nonumber\\
&\leq&  \frac{n_+^\beta}{n^\beta}\left(\sup_{t\in I_n}\left\{\frac{W(t)}{t}-\frac{W(n_+^\beta)}{n_+^\beta}\right\}
+\left|W(n_+^\beta)\right|\sup_{t\in I_n}
\left\{\left| \frac{1}{n_+^\beta}-\frac{1}{t}\right|\right\}
\right)
\rw 0,\: \:\mbox{a.s.}, \nonumber
\EQN
as $n\rw\IF$.
This confirms (\ref{e1}) and thus the proof has been completed. \endproof

\subsection{Proof of Theorem \ref{th.secUnstab}}
{\it Upper bound}:
Observe that
\begin{equation}
\label{upper1}
\pr(Q_1>u) \leq \pr(Q_1^{\phi_1 c} >u) \sim  \frac{1}{(\phi_1 c-\mu_1)(\alpha_1-1)}u^{1-\alpha_1}L_1(u),
\end{equation}
by Lemma \ref{th.Port}.

\vb

{\it Lower bound}:
Since
\begin{equation*}
\begin{split}
Q_1
& = \sup_{t \geq 0}\left\{Z_1(-t,0) + Z_2(-t,0) - ct - \sup_{s\in [0,t)}\{Z_2(-s,0) - C_2(-s,0)\} \right\} \\
& \geq \sup_{t \geq 0}\left\{Z_1(-t,0) + Z_2(-t,0) - ct - \sup_{s\in [0,t)}\{Z_2(-s,0) - \phi_2cs\} \right\},
\end{split}
\end{equation*}
for any $\varepsilon\in(0,1)$, application of Lemmas \ref{lemma2}-\ref{Lemma1} yields,
with $\lambda=\phi_2c$,
\begin{eqnarray}\label{star}
\lefteqn{
\pr(Q_1 > u)
\geq
\pr \left( \sup_{t \geq 0} \left\{ Z_1(-t,0) + Z_2(-t,0) - ct - \sup_{s\in [0,t)}\{Z_2(-s,0) - \phi_2cs\} \right\} > u \right)}\nonumber \\
&=&
\pr \left( \sup_{t \geq 0} \left\{ Z_1(t) - \phi_1 c t - \sup_{s \in [0,t)} \{ Z_2(s)-Z_2(t) - \phi_2c(s-t) \} \right\} > u \right) \nonumber\\
&\geq &
\pr \left( \sup_{0\leq t \leq {T_\lambda}(\varepsilon u)} \left\{ Z_1(t) - \phi_1 c t - \sup_{0\leq s \leq {T_\lambda}(\varepsilon u), v> s} \{ Z_2(s)-Z_2(v) - \phi_2c(s-v) \} \right\} > u \right) \nonumber\\
&= &
\pr \left( \sup_{0\leq t \leq {T_\lambda}(\varepsilon u)} \left\{ Z_1(t) - \phi_1 c t \right\} - \sup_{0\leq s \leq {T_\lambda}(\varepsilon u)} \check Q_2^{\phi_2c}(s) > u
\right)
\nonumber\\
&\geq&
\pr \left( \sup_{0\leq t \leq {T_\lambda}(\varepsilon u)} \left\{ Z_1(t) - \phi_1 c t  \right\} -  \sup_{0\leq s \leq {T_\lambda}(\varepsilon u)} \check Q_2^{\phi_2c}(s)> u ,
 \sup_{0\leq s \leq {T_\lambda}(\varepsilon u)} \check Q_2^{\phi_2c}(s)\leq \varepsilon u\right)\nonumber\\
&\geq&
\pr \left( \sup_{0\leq t \leq {T_\lambda}(\varepsilon u)} \left\{ Z_1(t) - \phi_1 c t \right\} > (1+\varepsilon)u \right)\pr\left( \sup_{0\leq s \leq {T_\lambda}(\varepsilon u)} \check Q_2^{\phi_2c}(s)\leq \varepsilon u\right)\nonumber\\
&\sim& \pr \left( Q_1^{\phi_1 c} > u \right), \ \ u\to\infty, \,\varepsilon\downarrow 0,\nonumber
\end{eqnarray}
which combined with Lemma \ref{th.Port}
leads to the asymptotic upper bound that matches the lower bound.
This completes the proof.
\halmos

\subsection{Proof of Theorem  \ref{th.secStabH}}
{\it Upper bound}:
Combining
\begin{equation}
\label{upper2}
\pr(Q_1>u) \leq \pr(Q_1 + Q_2 >u) = \pr \left( \sup_{t \geq 0} \{ Z_1(-t,0) + Z_2(-t,0) - ct \} > u \right)
\end{equation}
with $\pr(Z_1(1)+Z_2(1)>u)\sim \pr(Z_1(1)>u)$ as $u\to\infty$,
together with Lemma \ref{th.Port}, straightforwardly
gives that
\BQN\label{upper22}
\pr(Q_1>u)
\sim
\frac{1}{c-\mu}\int_u^{\infty}\pr(Z_1(1)>x){\rm d}x \sim \frac{1}{(c-\mu)(\alpha_1-1)}u^{1-\alpha_1}L_1(u),
\EQN
as $u\to\infty$.

\vb

{\it Lower bound}:
Let $\varepsilon >0$ be given.
Following the same argument as in the lower bound  of the proof of
Theorem \ref{th.secUnstab}, we have, with $\lambda=\mu_2^\varepsilon:= \mu_2-\varepsilon$,
\BQNY
&&\pr(Q_1 > u)\\
&& \ \ \geq\:
\pr \left( \sup_{t \geq 0}\left\{ Z_1(-t,0) + Z_2(-t,0) - ct - \sup_{s\in [0,t)} \left\{ Z_2(-s,0) - \phi_2cs \right\} \right\} > u \right) \\
&& \ \ =\:
\pr \left( \sup_{t\geq 0} \left\{ Z_1(t) - (c-\mu_2^\varepsilon)t -
           \sup_{s \in [0,t)} \{ Z_2(s)-Z_2(t) - \phi_2cs +\mu_2^\varepsilon t\} \right\} > u \right)\\
&&\ \ \geq\: \pr \left( \sup_{t\in [0,{T_\lambda}(\varepsilon u)]} \left\{ Z_1(t) - (c-\mu_2^\varepsilon)t -
           \sup_{s \in [0,{T_\lambda}(\varepsilon u)), s< t} \{ Z_2(s)-Z_2(t) +\mu_2^\varepsilon(t-s)\} \right\} > u \right)\\
&& \ \ \geq  \:\pr \left( \sup_{t\in [0,{T_\lambda}(\varepsilon u)]} \left\{ Z_1(t) - (c-\mu_2^\varepsilon)t -
           \sup_{s \in [0,{T_\lambda}(\varepsilon u)]}\check Q_2^{\mu_2^\varepsilon} (s)\right\} > u \right)\\
&& \ \ \geq \: \pr \left( \sup_{t\in [0,{T_\lambda}(\varepsilon u)]} \left\{ Z_1(t) - (c-\mu_2^\varepsilon)t
           \right\} > (1+\varepsilon)u \right)\pr\left(\sup_{s \in [0,{T_\lambda}(\varepsilon u)]}\check Q_2^{\mu_2^\varepsilon}(s) <\varepsilon u\right).
\EQNY
By Lemmas \ref{lemma2}--\ref{Lemma1} in combination with Lemma \ref{th.Port}, we obtain
$$\pr(Q_1 > u)\geq \pr(Q_1^{c-\mu_2^\varepsilon}>(1+\varepsilon)u)(1+o(1))\sim \frac{(1+\varepsilon)^{1-\alpha_1}}{(c-\mu+\varepsilon)(\alpha_1-1)}u^{1-\alpha_1}L_1(u), \ \ u\rw\IF.$$
Letting $\varepsilon\downarrow 0$, and recalling  (\ref{upper22}), completes the proof.
\halmos

\subsection{Proof of Theorem \ref{th.secStabL}}
{\it Upper bound}:
The starting point is the following evident equality:
$$Q_1 =\sup_{t \in [0, u^{1-\varepsilon})}U_1(t)\vee\sup_{t \geq u^{1-\varepsilon}}U_1(t),\:\:\:\:\:U_1(t):=
Z_1(-t,0) - C_1(-t,0).
$$
with $\varepsilon$ strictly between ${\alpha_1}/({1+\alpha_1})$ and $1$. Then, for
$$\mathscr{U}_{\varepsilon} \eq \{ \forall t\geq u^{1-\varepsilon} \,:\, Z_2(-t,0)+Q_2^{\phi_2c}(-t) \leq (\mu_2+\varepsilon)t \}$$
we have
\begin{equation}
\label{upper3a}
\pr \left( \sup_{t \geq u^{1-\varepsilon}}U_1(t)> u \right )=
\pr \left (  \sup_{t \geq u^{1-\varepsilon}}U_1(t) > u; \, \mathscr{U}_{\varepsilon}  \right )
 + \pr \left( \sup_{t \geq u^{1-\varepsilon}}U_1(t) > u; \, \mathscr{U}_{\varepsilon}^c  \right).
\end{equation}
It follows from (\ref{e2}) that on the event $\mathscr{U}_{\varepsilon}$, for $t\geq u^{1-\varepsilon}$,
\BQNY
B_2(-t,0)=Z_2(-t,0)+Q_2(-t)-Q_2(0)\leq Z_2(-t,0)+Q_2^{\phi_2c}(-t)\leq (\mu_2+\varepsilon)t,
\EQNY
 which together with the fact that $C_1(s,t)+B_2(s,t)=c(t-s)$ for all $s\leq t$ yields, for $t\geq u^{1-\varepsilon}$,
 $$C_1(-t,0)\geq (c-\mu_2-\varepsilon)t.$$
Moreover,
\begin{eqnarray}
\nonumber\label{upper3aa}
\pr \left( \sup_{t \geq u^{1-\varepsilon}}U_1(t) > u; \, \mathscr{U}_{\varepsilon} \right)
&\leq& \pr \left( \sup_{t\geq u^{1-\varepsilon}}\{Z_1(-t,0) -(c-\mu_2-\varepsilon)t \}>u;\, \mathscr{U}_{\varepsilon}  \right) \\
&\leq &\pr \left( \sup_{t\geq 0}\{Z_1(-t,0) -(c-\mu_2-\varepsilon)t\}>u \right) \pr( \mathscr{U}_{\varepsilon} );
\label{upper3aa}
\end{eqnarray}
the first term in (\ref{upper3aa}) is roughly of the order $u^{1-\alpha_1}$, whereas $\pr( \mathscr{U}_{\varepsilon}) \to  1$, as $u\to\infty$,
 as a consequence of the law of large numbers
and Lemma \ref{BIG}.

In addition,
\begin{eqnarray}\nonumber
\pr \left( \sup_{t \geq u^{1-\varepsilon}}U_1(t)> u; \, \mathscr{U}_{\varepsilon}^c  \right)
& \leq& \pr \left( \sup_{t\geq u^{1-\varepsilon}}\{Z_1(-t,0) -\phi_1ct \}>u;\, Y\mathscr{U}_{\varepsilon}^c  \right)\\
 &\leq  &\pr \left( \sup_{t\geq 0}\{Z_1(-t,0) -\phi_1ct\}>u \right) \pr(\mathscr{U}_{\varepsilon}^c ),
\label{upper3ab}
\end{eqnarray}
where the first term in (\ref{upper3ab}) essentially vanishes as $u^{1-\alpha_1}$, but $\pr(\mathscr{U}_{\varepsilon}^c ) \to 0$ due to the law of large numbers
and Lemma \ref{BIG}. We conclude it is negligible relative to (\ref{upper3aa}).

Combining (\ref{upper3aa}) and (\ref{upper3ab}) gives that
\[
\pr \left(\sup_{t \geq u^{1-\varepsilon}}\{ Z_1(-t,0) - C_1(-t,0)\} > u \right) \leq \pr(Q_1^{c-\mu_2-\varepsilon}>u ), \quad \textrm{as} \, u\to \infty.
\]
Now we are left with showing that
\[
\pr(Q_1 > u) \sim \pr \left (\sup_{t \geq u^{1-\varepsilon}}U_1(t)>u \right), \,\:\: \textrm{as} \:u \to \infty.
\]
Since we have that
\begin{equation}
\label{negl}
\pr(Q_1>u) = \pr \left( \sup_{t\geq 0}U_1(t)>u \right)  \leq
\pr \left( \sup_{t \geq u^{1-\varepsilon}}U_1(t)>u \right)
+ \pr \left(\sup_{t\in [0, u^{1-\varepsilon}]}U_1(t)>u \right)
\end{equation}
and we have already showed that both
$\pr\left( \sup_{t \geq u^{1-\varepsilon}}U_1(t)>u \right)$ and $\pr(Q_1 >u)$  are of order $u^{1-\alpha_1}$,
we see that it suffices to prove that the last term in (\ref{negl}) is negligible.

Let
$$S_n:=\sup_{s\in [n,n+1]}Z_1(s)-Z_1(n)-\phi_1c(s-n), \:\:\:n\in \mathbb{N}.$$
Using that $S_n, n\in \mathbb{N}$ are i.i.d., we get
\BQNY
\pr \left( \sup_{t\in [0, u^{1-\varepsilon}]}U_1(t)>u \right)
&\leq& \pr \left( \sup_{t\in [0, u^{1-\varepsilon}]}\{ Z_1(-t,0)-\phi_1ct)\}>u \right) \\
& = &\pr \left( \sup_{t\in [0,u^{1-\varepsilon}]}\{Z_1(t)-\phi_1ct)\}>u \right)\\
& \leq &\pr \left( \sum_{i=0}^{[u^{1-\varepsilon}]} S_i>u \right)\leq ([u^{1-\varepsilon}]+1)\pr\left(S_0>\frac{u}{[u^{1-\varepsilon}]+1}\right).
\EQNY
Hence, Lemma \ref{th.Willekens}, in combination with the fact that $\varepsilon$ lies strictly between ${\alpha_1}/({1+\alpha_1})$ and $1$, leads to
\BQNY
\pr \left( \sup_{s\in [0, u^{1-\varepsilon}]}U_1(s)>u \right)\leq \bar \kappa \,L_1(u^{\varepsilon})u^{1-\varepsilon-\varepsilon \alpha_1}=o\left(L_1(u)u^{1-\alpha_1}\right),
\EQNY
where $\bar\kappa$ is a positive constant.
This completes the proof of the upper bound.

\vb

{\it Lower bound}:
The proof of the lower bound is the same as in the proof of Theorem \ref{th.secStabH}. 
Relying on Lemma \ref{th.Port}, we then obtain the equivalence of the asymptotic upper and lower bound.
\halmos

\subsection{Proof of Theorem  \ref{th.reqc}  }
Since the proof of this scenario needs a case-specific approach that involves the notion of tandem systems,
we begin with some notation and auxiliary results.

For $\varepsilon$ such that $\phi_1 c-\mu_1<\varepsilon<c-\mu$, let
$$V^{\varepsilon} \eq \sup_{t \geq 0}\{Z_2(-t,0) - (c-\mu_1 - \varepsilon)t\} - \sup_{s \geq 0} \{Z_2(-s,0) - \phi_2cs\}.$$
Recall that
$$Q_1^d \eq \sup_{ t \geq 0}\{Z_1(-t,0) - dt\}, \quad d>\mu_1$$
and introduce
$$ \check Q_1^d \eq \sup_{t \geq 0}\{dt - Z_1(-t,0)\}, \quad d<\mu_1.$$

The following lemma states a straightforward counterpart of Lemma 2.1 in  \cite{DMGPS}.
\begin{Lemma}
\label{lemtand}
For $\varepsilon >0$ small enough, any $u$ and $x$, and for $\delta \in (0,1)$, we have
$$\pr(V^{-\varepsilon} > u+x)\pr(\check Q_1^{\mu_1 - \varepsilon}\leq x) \leq \pr(Q_1>u) \leq \pr(V^{\varepsilon}>(1-\delta)u) + \pr(Q_1^{\mu_1 +\varepsilon}>\delta u)$$
\end{Lemma}

A combination of  Theorem 4.7 in  \cite{LiM} (see also Theorem 12.9 in \cite{DeM15})
with Lemma \ref{th.Port} leads to the following lemma.
\begin{Lemma}
\label{tandAs}
Let  $|\varepsilon|<\min(c-\mu,\mu_1- \phi_1 c)$ and $Z_2$ be spectrally positive
with $\alpha_2\notin \N$. Then, as $u\to\infty$,
\[
\pr(V^\varepsilon>u)\sim
\left( \frac{\mu_1-\phi_1c +\varepsilon}{\phi_2c-\mu_2} \right)^{\alpha_2-1}\frac{1}{(c-\mu-\varepsilon)(\alpha_2-1)}u^{1-\alpha_2}L_2(u) .
\]
\end{Lemma}

\textit{Proof of Theorem \ref{th.reqc}}:
Let $\delta\in (0,1)$ and $\varepsilon>0$ be such that $\varepsilon<\min(c-\mu,\mu_1- \phi_1 c)$.
Then, following Lemma \ref{tandAs}, as $u\to\infty$,
\begin{eqnarray*}
\pr(\check Q_1^{\mu_1-\varepsilon}<\sqrt{u})&\to& 1,\\
\pr(V^{-\varepsilon}>u+\sqrt{u})
&\sim&
\left( \frac{\mu_1-\phi_1c -\varepsilon}{\phi_2c-\mu_2} \right)^{\alpha_2-1}\frac{1}{(c-\mu+\varepsilon)(\alpha_2-1)}u^{1-\alpha_2}L_2(u)\\
\pr(V^{\varepsilon}>(1-\delta)u)
&\sim&
\left( \frac{\mu_1-\phi_1c +\varepsilon}{\phi_2c-\mu_2} \right)^{\alpha_2-1}\frac{1}{(c-\mu-\varepsilon)(\alpha_2-1)}(1-\delta)^{1-\alpha_2}u^{1-\alpha_2}L_2(u).
\end{eqnarray*}
Since $\alpha_1>\alpha_2$, we find by Lemma \ref{th.Port} that for each $\delta \in(0,1)$,
\begin{eqnarray*}
\pr(Q_1^{\mu_1+\varepsilon}>\delta u)=o(\pr(V^{\varepsilon}>(1-\delta)u)),
\end{eqnarray*}
as $u\to\infty$.
Thus, by Lemma \ref{lemtand}, passing with $\varepsilon,\delta\downarrow 0$,
we obtain that
\[
\pr(Q_{1} > u) \sim \left( \frac{\mu_1-\phi_1c}{\phi_2c-\mu_2} \right)^{\alpha_2-1}\frac{1}{(c-\mu)(\alpha_2-1)}u^{1-\alpha_2}L_2(u).
\]
This completes the proof.
\halmos

\section*{Acknowledgments}
{\small
K. D\c ebicki
was partially supported by NCN Grant No 2015/17/B/ST1/01102 (2016-2019) whereas
P. Liu was partially supported by  the Swiss National Science Foundation Grant 200021-166274.
{M.\ Mandjes'} research is partly funded by the NWO Gravitation project N{\sc etworks}, grant number 024.002.003. He is also affiliated to (A)~CWI, Amsterdam, the Netherlands; (B)~E{\sc urandom}, Eindhoven University of Technology, Eindhoven, the Netherlands;
and (C)~Amsterdam Business School, Faculty of Economics and Business, University of Amsterdam,
Amsterdam, the Netherlands.}

\end{document}